# Packing Ovals In Optimized Regular Polygons

Frank J. Kampas [a], János D. Pintér [b], Ignacio Castillo [c, *]

[a] Physicist at Large Consulting LLC, Bryn Mawr, PA, USA
[b] Department of Industrial and Systems Engineering, Lehigh University, Bethlehem, PA, USA
[c] Lazaridis School of Business and Economics, Wilfrid Laurier University, Waterloo, ON, Canada
[*] Corresponding author. E-mail address: icastillo@wlu.ca

**Abstract**

We present a model development framework and numerical solution approach to the general problem-class of packing convex objects into optimized convex containers. Specifically, here we discuss the problem of packing ovals (*egg-shaped* objects, defined here as generalized ellipses) into optimized regular polygons in $\mathbb{R}^2$. Our solution strategy is based on the use of embedded Lagrange multipliers, followed by nonlinear (global-local) optimization. The numerical results are attained using randomized starting solutions refined by a single call to a local optimization solver. We obtain credible, tight packings for packing 4 to 10 ovals into regular polygons with 3 to 10 sides in all (224) test problems presented here, and for other similarly difficult packing problems.

## 1 Introduction

The efficient packing, arrangement, or configuration design of objects is required across a vast range of engineering and scientific applications. Industrial engineering (IE) and operations research (OR) applications include facility layout design, cutting stock problems in various contexts (e.g., in the glass, metal, paper, textile, and wood industries), and the packing arrangement of solid objects for storage or transportation. A perfunctory web search for the key words "packing problems" returns more than 150,000 results (as of November 2018). Therefore here we only refer to a relatively small selection of works that present, review and discuss IE/OR packing methods and applications: consult e.g., Dowsland and Dowsland (1992), Lodi *et al.* (2002), Pisinger and Sigurd (2007), Bennell and Oliveira (2008), Castillo *et al.* (2008), Hifi and M'Hallah (2009), Bennell *et al.* (2010), Chernov *et al.* (2010), López and Beasley (2011), Fasano (2014), Fasano and Pintér (2015), Alt (2016), Anjos and Vieira (2017), with many topical references. Let us add that e.g., information theory and error-correcting codes (Conway, 1995; Shannon, 1948), number theory, approximation theory, algebra, and theoretical physics (Cohn, 2010), and the design of experiments (Kleijnen, 2015) are further important application areas.

As stated by Saunders (2017) in his anniversary review article, Thompson's classical work (1917) has stood at the forefront of our understanding of the development of biological form for the past century. Studies in computational physics, chemistry and biology often search for optimized object configurations. The number of web references related to these topics is in the range of *millions*: here we refer only to a small selection of books by Jensen (2007), Landau *et al.* (2012), Newman (2012), Kellis (2016), and O'Neil (2017).

To give a few concrete examples, object packing studies in the sciences are related, e.g., to the structural analysis of liquids, crystals, and glasses (Bernal, 1959); the flow and compression of granular materials (Edwards, 1994; Jaeger and Nagel, 1992, 1996); the design of high-density ceramic materials and the formation and growth of crystals (Cheng *et al.*, 1999; Rintoul and Torquato, 1996); the thermodynamics of liquid to crystal transition (Alder and Wainwright, 1957; Chaikin, 2000; Pusey, 1991); and the chromosome organization in human cell nuclei (Uhler and Wright, 2013).

The paradigm packing objects *efficiently* leads to interesting model development and (often) hard optimization challenges. A packing configuration can be formally defined as a non-overlapping arrangement of a given collection of objects inside a chosen type(s) of container(s). Packings can be optimized according to some appropriately selected criterion, such as the area or volume of the container or the packing fraction (defined as the fraction of the container area/volume covered by the packed objects). The convexity of the packed objects and/or of the container is often postulated, noting that such problems typically still require non-convex continuous and/or combinatorial optimization approaches.

To start with a seemingly "easy" case, packing identical circles has received considerable attention. Research on packing identical circles frequently aims at proving the optimality of the configurations found, either theoretically or with the help of rigorous computational techniques. Provably optimal configurations, with the exception of certain special cases, are available only for a few tens of circles; best-known results are available for packing up to 2,600 identical circles in a circle and 10,000 identical circles in a square. For further details and references, consult e.g., Szabó *et al.* (2007) or Specht (2018).

The general circle packing problem – considered for a given set of circles with (in principle) arbitrary size – is a substantial generalization of the case with identical circles. In full generality, provably optimal configurations are available only for models with $n \leq 4$ circles. Therefore, studies dealing with such problems introduce and apply efficient generic or tailored global scope numerical solution strategies, but without the ability to prove the optimality of the results obtained. We refer to Castillo *et al.* (2008) and to Hifi and M'Hallah (2009) for reviews of both uniform and arbitrary sized circle packings and their applications. More recently, Pintér *et al.* (2017) present numerical results for general sphere packings in 2, 3, 4, 5 dimensions with up to 50 spheres.

Compared to circle packings, ellipse packing problems have received relatively little attention in the literature. Finding high quality, globally optimized packings of ellipses with arbitrary size and orientation is a hard computational problem. The key challenge is the modeling and enforcement of the non-overlapping constraints as a function of ellipse center locations and orientations. Galiev and Lisafina (2013) studied of ellipse packing problems assuming identical ellipses, orthogonally oriented inside a rectangular container. Uhler and Wright (2013) relax these assumptions, and propose a model that minimizes a measure of overlap between ellipses (while overlaps still remain possible). Kallrath (2017) extends the work presented in Kallrath and Rebennack (2014) to pack non-overlapping ellipses of arbitrary size and orientation into optimized rectangular containers: the key modeling idea is to use separating lines to ensure that the ellipses do not overlap with each other. Birgin *et al.* (2017) extend the work presented in Birgin *et al.* (2016) for packing arbitrary ellipses in convex containers: they propose a multi-start strategy combined with starting guesses and a local optimization solver, in order to find good quality packings with up to 1000 ellipsoids. We will refer to our related previous work in Section 3, while extending that work for the more general problem-class studied here. Despite the substantial amount of research effort highlighted above, it is evident that packing relatively "simple" objects such as circles or ellipses into "simple" (circular, rectangular) containers already leads to modeling and optimization challenges.

In this article, we concentrate on model development and numerical solution approaches to the problem of packing more general convex objects using a unified, flexible, accurate, and efficient framework. Our model development relies on representing the objects to be packed by a general convex set-type known as ovals or eggs. (In common English, the term oval is used for any shape which reminds one of an egg.) We will give a formal definition of eggs, noting that – according to our definition – circles and ellipses are special cases of eggs. To further enhance the flexibility of our modeling framework, we consider optimized regular polygons as container sets.

To our best knowledge, there are no previous studies related to the general problem-class studied here. Following this introduction, we define and discuss eggs in Section 2. The optimization model development framework is presented in Section 3. Illustrative numerical results and their analysis are described in Section 4. The conclusions (Section 5) are followed by a fairly extensive list of references.

## 2 Eggs: General Definition and Some Special Cases

Perhaps surprisingly, there is no unique definition of ovals, egg shaped curves or bodies. Consult e.g. the topical webpages of Köller (2018) and Yamamoto (2018) for related discussions and a number of alternative definitions.

Citing the "oval" entry of Wikipedia (2018), oval commonly means a shape like an egg or an ellipse. It can be also used to refer to a "stadium" shape defined by two semicircles joined by a rectangle. Sometimes, it can even refer to a rectangle with rounded corners. The general term oval is used also to describe certain non-convex objects such as the lemniscate of Bernoulli, or Cassini ovals.

In this work, we define the contour of an egg curve in $\mathbb{R}^2$ that follows the equation

$$\left(\frac{x}{a}\right)^p + e^{tx}\left(\frac{y}{b}\right)^p = 1. \qquad (1)$$

In (1), $(x, y)$ is the location of a point on the egg contour curve in $\mathbb{R}^2$, $a > 0$ and $b > 0$ are the semi-major and semi-minor axes of this curve, $p \geq 2$ is an even integer, and $t \geq 0$ is the distortion factor. An egg can be defined with arbitrary size and orientation parameters, noting that – in the case studied here –the size and distortion factors are limited in order to maintain the egg's convexity. The four images displayed in Figure 1 illustrate the modeling flexibility offered by eggs defined in the form of equation (1), for several settings of the parameters $a$, $b$, $p$, and $t$. For example, the setting $a = b$, $p = 2$, $t = 0$ results in a circle; and the (general) pair of $a$, $b$ with $p = 2$, $t = 0$ results in an ellipse. Setting $p = 2$ leads to ovals which become more asymmetrical and "pointed" as $t$ increases. Setting $p = 4$ leads to "stadium" like objects, and so on.

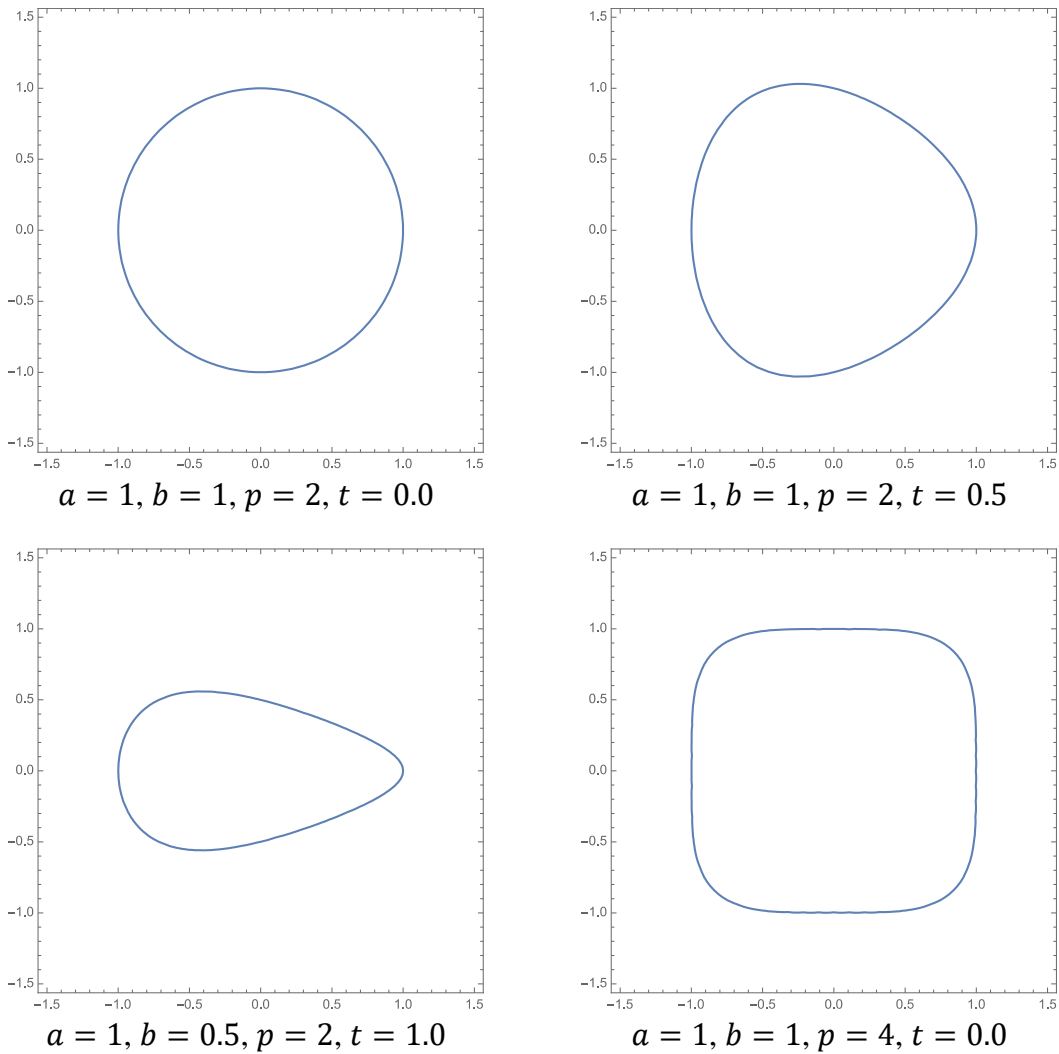

Figure 1. General egg curve examples for several values of $a$, $b$, $p$, and $t$.

## 3 Optimization Model Development

Our modeling and solution approach is based on extending the use of embedded Lagrange multipliers to the case of egg objects studied here. Embedded Lagrange multipliers were introduced to pack ellipses inside optimized circular containers in Kampas *et al.* (2017), and to pack ellipses inside optimized regular polygons in Kampas *et al.* (2018).

Here we also use an optimized regular polygon as the container set. The extension from ellipses to eggs is, however, not trivial. Our objective is to minimize the area of the regular polygon that contains a given collection of eggs with arbitrary size and orientation.

The input data to such an optimization problem instance are defined by the number of sides for the container, and the parameters (semi-major, semi-minor axes, exponent and distortion factor) of the eggs to be packed.

The primary decision variables are the polygon's apothem, and the centre position and orientation of each of the packed eggs. (Recall that the apothem of a regular polygon is a line segment from the polygon's center to the midpoint of one of its sides.)

There are two sets of secondary variables. The first set consists of the positions of the distance maximizing lines pointing from each egg boundary to the center of each of the polygon faces. The second set is given by the positions of the points on one of each pair of eggs which minimizes the value of the equation describing the other egg.

The secondary variables are used to define the model constraints. The first set is used to represent the constraints that keep the eggs inside the container. The second set is used to prevent the eggs from overlapping. These constraint sets are generated by embedded Lagrange multiplier conditions. In our optimization strategy, the calculations for finding the polygon's apothem and for preventing egg overlaps proceed simultaneously, rather than being performed to completion at each step towards the minimization of the polygon area.

Denote by $m$ the number of polygon sides, and by $d$ the $m$-polygon's apothem: the area of this polygon equals $m \cdot d^2 \cdot \tan(\pi/m)$. Since $m$ is an input parameter to the optimization problem instance, minimizing the apothem is equivalent to minimizing the area of the regular $m$-polygon.

Equation $e(a, b, p, t, xc, yc, \theta; x, y)$ displayed below defines an egg with its given input parameters – semi-major and semi-minor axes $a$ and $b$, exponent $p$, and distortion factor $t$ – centered at $\{xc, yc\}$, and rotated counterclockwise by angle $\theta$. Recall that $(xc, yc)$ and $\theta$ are the primary decision variables for each egg.

$$e(a, b, p, t, xc, yc, \theta; x, y) \tag{2}$$
$$= \left(\frac{\cos(\theta)(x - xc) + \sin(\theta)(y - yc)}{a}\right)^p$$
$$+ e^{t\delta}\left(\frac{\cos(\theta)(y - yc) - \sin(\theta)(x - xc)}{b}\right)^p - 1,$$
$$\text{where } \delta = \cos(\theta)(x - xc) + \sin(\theta)(y - yc).$$

The value of $e(a, b, p, t, xc, yc, \theta; x, y)$ is negative for all points $(x, y)$ located inside the egg, zero for all points on the egg boundary, and positive for all points outside the egg.

We will assume (postulate) that the optimized polygon container is centered at the origin. Consider now the line $cx \cdot x + cy \cdot y = l$ that embeds one of the sides of the polygon. The slope of this line is $-(cx/cy)$. If $(cx, cy)$ is a unit vector so that $cx^2 + cy^2 = 1$, then the point on the line closest to the origin is $l \cdot (cx, cy)$. The slope of the line to that point is $cy/cx$: hence, the line from the origin to the closest point on $cx \cdot x + cy \cdot y = l$ is perpendicular to it.

As mentioned above, the first set of constraints is used to keep the eggs inside the container. We find the maximum value of $l$, denoted by $l^{\max}$, for which the side of the polygon intersects the egg. Our $m$-polygon derivation follows the first order Karush–Kuhn–Tucker conditions described by Kallrath and Rebennack (2014) for a rectangular container. For the egg to be contained inside the polygon, all sides $cx \cdot x + cy \cdot y$ must be less than or equal to this maximum value. Consider the following equation using $\lambda$ as the maximizing embedded Lagrange multiplier.

$$cx \cdot x + cy \cdot y = \lambda \cdot e(a, b, p, t, xc, yc, \theta; x, y). \tag{3}$$

Differentiating equation (3) with respect to $x$, $y$, and $\lambda$, we can numerically obtain the maximum value $l^{\max}$ in the direction $(cx, cy)$.

For example, consider the line defined by $x + y = 1$ and the egg defined by $e(1/2, 1/3, 2, 0, 0, 0, 1/2; x, y)$. The maximum value in the direction $(cx, cy)$ is $l^{\max} = 0.2444$ with $(x, y) = (0.4000, 0.2914)$: this value can be found by numerically solving equation (3). Figure 2 shows the corresponding line-egg configuration.

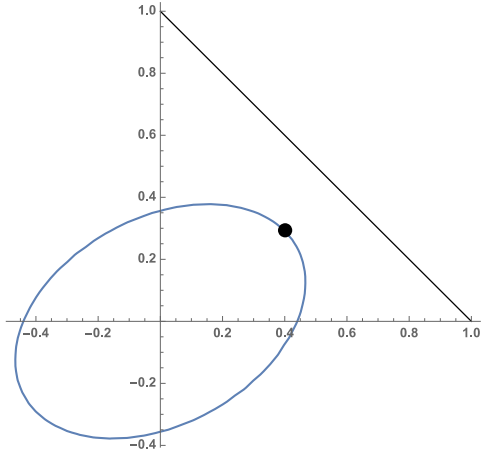

Figure 2. Line-egg configuration in the given example.

Based on the above, the condition of containing a given egg inside the polygon can be described by the relation

$$d \geq l^{\max}(a,b,p,t,xc,yc,\theta;x,y). \tag{4}$$

This constraint will be used in our optimization framework, for all sides of the regular polygon which, of course, share the same apothem distance $d$. Thus, for a regular polygon with $m$ sides, the points $(cx_k, cy_k)$ that define the unit vectors for each apothem are given by

$$(cx_k, cy_k) = \left(\cos\left(\frac{2 \cdot k \cdot \pi}{m} - \frac{\pi}{2}\right), \sin\left(\frac{2 \cdot k \cdot \pi}{m} - \frac{\pi}{2}\right)\right), k = 1, \ldots, m. \tag{5}$$

Proceeding next towards preventing egg overlaps, it is useful to determine equations for the derivatives of the egg equation with respect to $x$ and $y$. For a given set of values $(a,b,p,t,xc,yc,\theta)$ we can denote $e(a,b,p,t,xc,yc,\theta;x,y)$ simply as $e(x,y)$: then these derivatives are

$$\frac{de(x,y)}{dx} = -\frac{e^{t\delta}p\varphi^{p-1}}{b}\sin(\theta) + e^{t\delta}t\varphi^p \cos(\theta) + \frac{p\omega^{p-1}}{a}\cos(\theta), \tag{6}$$

$$\frac{de(x,y)}{dy} = \frac{e^{t\delta}p\varphi^{p-1}}{b}\cos(\theta) + e^{t\delta}t\varphi^p \sin(\theta) + \frac{p\omega^{p-1}}{a}, \tag{7}$$

where

$$\varphi = \left(\frac{1}{b}\right)\left((y - yc)\cos(\theta) + (xc - x)\sin(\theta)\right), \tag{8}$$

$$\omega = \left(\frac{1}{a}\right)\left((x - xc)\cos(\theta) + (y - yc)\sin(\theta)\right). \tag{9}$$

All pairs of packed eggs are prevented from overlapping by requiring that the minimum value of the egg equation for egg $i$ for any point on egg $j$ has to be greater than a sufficiently small parameter $\varepsilon \geq 0$. This requirement will be met by using embedded Lagrange multiplier conditions.

The equations shown below determine the point on egg $j$ that maximizes or minimizes the value of the function describing egg $i$. In the case considered here, $\lambda$ must be negative to obtain the minimum. During optimization, this requirement with respect to the sign of $\lambda$ will be enforced by setting its search bounds.

$$\begin{cases} \frac{de(x,y)}{dx}(a_i, b_i, p_i, t_i, xc_i, yc_i, \theta_i; x, y) = \lambda \cdot \frac{de(x,y)}{dx}(a_j, b_j, p_j, t_j, xc_j, yc_j, \theta_j; x, y) \\ \frac{de(x,y)}{dy}(a_i, b_i, p_i, t_i, xc_i, yc_i, \theta_i; x, y) = \lambda \cdot \frac{de(x,y)}{dy}(a_j, b_j, p_j, t_j, xc_j, yc_j, \theta_j; x, y) \\ e(a_j, b_j, p_j, t_j, xc_j, yc_j, \theta_j)(x, y) = 0 \end{cases} \quad (10)$$

The last equation type introduced corresponds to the requirement that the distance minimizing point lies on egg $j$. Eliminating $\lambda$ from the first two equations, we obtain

$$\begin{aligned} \frac{de(x,y)}{dy}&(a_i, b_i, p_i, t_i, xc_i, yc_i, \theta_i; x, y) \\ &\cdot \frac{de(x,y)}{dx}(a_j, b_j, p_j, t_j, xc_j, yc_j, \theta_j; x, y) \\ &= \frac{de(x,y)}{dy}(a_j, b_j, p_j, t_j, xc_j, yc_j, \theta_j; x, y) \\ &\cdot \frac{de(x,y)}{dx}(a_i, b_i, p_i, t_i, xc_i, yc_i, \theta_i; x, y). \end{aligned} \quad (11)$$

Let us note that at the point on egg $j$ that minimizes or maximizes the value of the function describing egg $i$, the slope of egg curve $i$ equals the slope of egg curve $j$.

For example, consider two eggs defined by $e_i(1,1,2,1/2,0,0,0; x, y)$ and $e_j(1/2,3/4,2,0,0,0,0; x, y)$. The overlap indicator value between the eggs is $-0.5924$ with $(x, y) = (0.1035, 0.7338)$: this can be found by solving equations (10). Figure 3 shows this (overlapping eggs) configuration.

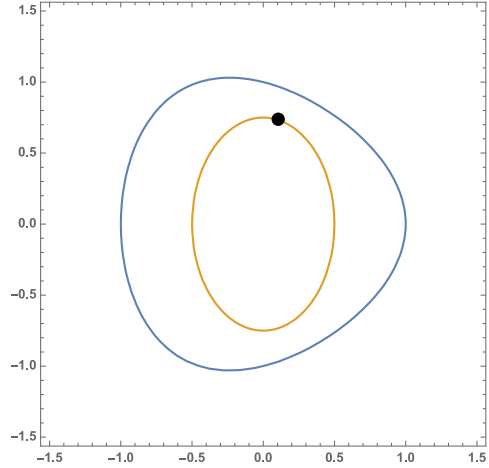
Figure 3. Two overlapping eggs.

To give another example, consider two eggs defined by $e_i(3/4,1,2,1/2,-1/2,0,0;x,y)$ and $e_j(1/3,1/2,2,0,3/4,0,1/2;x,y)$. The overlap indicator value between the eggs is 0.5696 with $(x,y) = (0.3745, 0.1004)$: again, this value can be found by solving the set of equations (10) for $x$ and $y$. Figure 4 shows the resulting (non-overlapping eggs) configuration.

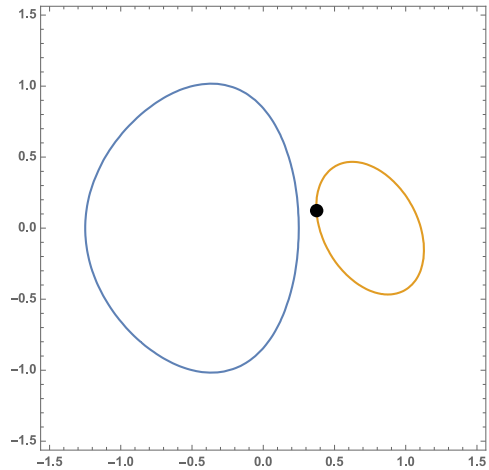
Figure 4. Two non-overlapping eggs.

In our optimization framework, $\lambda_{i,j}$ denote the Lagrange multipliers appearing in the equations to find the point $(x_{j,i}, y_{j,i})$ on egg $j$ that minimizes the value of the equation describing egg $i$. This calculation is restricted to minimization by requiring that the value of $\lambda_{i,j}$ is negative. Finally, we state constraints to prevent egg $i$ from overlapping with egg $j$, by requiring that the minimal value of the equation describing egg $i$ at the minimizing point on egg $j$ is at least $\varepsilon$.

Summarizing the model components and development steps discussed above, we present the following model-class for packing $n$ egg-shaped objects (ovals) in an optimized regular $m$-polygon.

$$\text{minimize} \quad d \tag{12}$$

$$\begin{aligned}
\text{subject to} \quad & d \geq l_{i,k}^{\max}(a_i, b_i, p_i, t_i, xc_i, yc_i, \theta_i; cx_k, cy_k) && \text{for } i = 1, \ldots, n \\
& && k = 1, \ldots, m \\[4pt]
& \frac{de(x,y)}{dx}(a_i, b_i, p_i, t_i, xc_i, yc_i, \theta_i; x_{j,i}, y_{j,i}) && \text{for } i = 1, \ldots, n-1 \\
& = \lambda_{j,i} && j = i+1, \ldots, n \\
& \cdot \frac{de(x,y)}{dx}(a_j, b_j, p_j, t_j, xc_j, yc_j, \theta_j; x_{j,i}, y_{j,i}) \\[4pt]
& \frac{de(x,y)}{dy}(a_i, b_i, p_i, t_i, xc_i, yc_i, \theta_i; x_{j,i}, y_{j,i}) && \text{for } i = 1, \ldots, n-1 \\
& = \lambda_{j,i} && j = i+1, \ldots, n \\
& \cdot \frac{de(x,y)}{dy}(a_j, b_j, p_j, t_j, xc_j, yc_j, \theta_j; x_{j,i}, y_{j,i}) \\[4pt]
& e(a_j, b_j, p_j, t_j, xc_j, yc_j, \theta_j; x_{j,i}, y_{j,i}) = 0 && \text{for } i = 1, \ldots, n-1 \\
& && j = i+1, \ldots, n \\[4pt]
& e(a_i, b_i, p_i, t_i, xc_i, yc_i, \theta_i) \geq \varepsilon && \text{for } i = 1, \ldots, n-1 \\
& && j = i+1, \ldots, n \\[4pt]
& lb \leq xc_i \leq ub && \text{for } i = 1, \ldots, n \\[4pt]
& lb \leq yc_i \leq ub && \text{for } i = 1, \ldots, n \\[4pt]
& -\pi \leq \theta_i \leq \pi && \text{for } i = 1, \ldots, n \\[4pt]
& lb \leq x_{j,i} \leq ub && \text{for } i = 1, \ldots, n-1 \\
& && j = i+1, \ldots, n \\[4pt]
& lb \leq y_{j,i} \leq ub && \text{for } i = 1, \ldots, n-1 \\
& && j = i+1, \ldots, n \\[4pt]
& 2 \cdot lb \leq \lambda_{j,i} \leq 0 && \text{for } i = 1, \ldots, n-1 \\
& && j = i+1, \ldots, n
\end{aligned}$$

In model (12) $lb$ and $ub$ are lower and upper bounds for the egg center positions: these bounds can be chosen appropriately for each egg packing model instance, in order to facilitate the finding of feasible solutions.

We close this section by noting that the embedded Lagrange multiplier conditions used require that the value of the function describing one egg is always greater than (or less than) 0 on all points on the other egg (or line), to ensure that eggs do not overlap with other eggs and that they are located within the polygon container. Analogous conditions are applicable to *any* convex curve described by a function which is always positive for points outside the curve: thus our approach can be applied also to packing rounded polygons.

## 4 Illustrative Numerical Results

We used a PC running under Windows 7, with an Intel Core i5 processor running at 2.6 GHz, with 16 GBytes of RAM, using *Mathematica* version 11. For the purposes of this study, we have implemented a rather efficient, but simple global-local optimization strategy based on multiple starting points. The numerical experiments are based on using a randomized starting solution used by a single call to the local nonlinear optimization solver Ipopt (2018). While – similarly to all other researchers who address similarly (or even less) difficult general packing problems – we cannot guarantee the theoretical optimality of the configurations found, our computational results consistently lead to visibly good quality packings.

As we stated earlier, to our best knowledge, there are no previously studied model instances available for the general egg packing problem-class considered here. Hence, we created our own test models, by packing eggs of (in principle) arbitrary size and orientation. Table 1 summarizes the egg packing test problem sets considered. Note that setting $p = 2$, $a_i = b_i$, and $t = 0.0$, test case 1 becomes a general circle packing problem for circles with radii $a_i = i^{-1/2}$. With $p = 2$, $a_i > b_i$, and $t = 0.0$, test case 2 becomes a general ellipse packing problem for ellipses with semi-major $a_i = i^{-1/2}$ and semi-minor $b_i = a_i/c$ axes, with eccentricity $c = 2$. Test cases 1-6 and 8 consider eggs with the same distortion factor set ($t = 0.0$, 0.5, or 1.0). Test case 7 considers eggs with different distortion factors $t_i = i/5$.

Tables 2-9 summarize the computational results for a total of 224 egg packing problem instances solved. In these model instances, the number of egg objects is chosen as one of $n = 4, \ldots, 10$; the number of regular polygon container sides is chosen as one of $m = 3, 4, 5,$ and 10.

Given how difficult these problems are thought to be, we believe that the CPU times (ranging from seconds to several hours, see result tables) are reasonable for our entire range of parameter choices. Note that the result tables also include container area and packing fraction information. Several illustrative packing configurations are displayed in Figure 5. (Similar figures are available, for all examples presented here.)

Considering Tables 2-7, we see that, on average, CPU time increases with the distortion factor $t$ (comparing, e.g., Tables 4-5 and 6-7). On average, CPU time also increases with the exponent $p$ (comparing Tables 2 and 8). We can also see that, on average, the packing fraction increases with the eccentricity $c$ (comparing, e.g., Tables 2 and 3). In general, there is a positive correlation between the packing fraction and certain functions of the different input parameters. Multiple linear regression analysis, using, for consistency, Tables 2-7

(with $p = 2$ and $t_i = t$ for all instances), indicates that the regression function for the packing fraction can be estimated as

$0.8022 + 0.0173 \cdot t + 0.0444 \cdot c - 0.2745 \cdot (1/n) - 0.1905 \cdot (1/m)$.

We received *p*-values (i.e., observed significance levels, cf. Black *et al.* (2014)) below 0.0001 for all other input parameters. This finding indicates that we have very strong statistical evidence suggesting that the regression coefficients are different from zero. Figure 6 highlights the quality of our regression equation by depicting actual vs. predicted packing fractions.

To conclude the discussion of numerical results, we note that it is apparent that our global-local optimization strategy could become prohibitively expensive for arbitrarily increasing size collections of eggs, particularly when the eggs are highly distorted and eccentric. Therefore we believe that large-scale problem instances open future research directions towards finding suitable heuristic approaches.

Table 1. Egg packing test cases

| Test case | $p_i$ | $(a_i, b_i)$ | $t_i$ |
|---|---|---|---|
| 1 | 2 | $(i^{-1/2}, a_i)$ | 0.0 |
| 2 | 2 | $(i^{-1/2}, a_i/2)$ | 0.0 |
| 3 | 2 | $(i^{-1/2}, a_i)$ | 0.5 |
| 4 | 2 | $(i^{-1/2}, a_i/2)$ | 0.5 |
| 5 | 2 | $(i^{-1/2}, a_i)$ | 1.0 |
| 6 | 2 | $(i^{-1/2}, a_i/2)$ | 1.0 |
| 7 | 2 | $(i^{-1/2}, a_i)$ | $i/5$ |
| 8 | 4 | $(i^{-1/2}, a_i)$ | 0.0 |

Table 2. Egg packing results for test case 1

| Problem number | Number of container sides | Number of eggs | Objective | Area of optimized container | Packing fraction | Time (sec.) |
|---|---|---|---|---|---|---|
| 1 | 3 | 4 | 2.6781 | 9.3169 | 0.7025 | 21.7489 |
| 2 | | 5 | 2.8097 | 10.2554 | 0.6995 | 15.5796 |
| 3 | | 6 | 2.8829 | 10.7962 | 0.7129 | 279.985 |
| 4 | | 7 | 2.9296 | 11.1491 | 0.7306 | 71.5815 |
| 5 | | 8 | 3.0304 | 11.9295 | 0.7157 | 589.264 |
| 6 | | 9 | 3.0490 | 12.0760 | 0.7360 | 67.2593 |
| 7 | | 10 | 3.0839 | 12.3544 | 0.7448 | 69.9138 |
| 8 | 4 | 4 | 2.1580 | 9.3137 | 0.7027 | 14.0704 |
| 9 | | 5 | 2.1871 | 9.5668 | 0.7498 | 54.0113 |
| 10 | | 6 | 2.2897 | 10.4850 | 0.7341 | 28.3510 |
| 11 | | 7 | 2.3275 | 10.8340 | 0.7519 | 62.3785 |
| 12 | | 8 | 2.3561 | 11.1021 | 0.7691 | 57.7027 |
| 13 | | 9 | 2.3702 | 11.2353 | 0.7910 | 134.182 |
| 14 | | 10 | 2.4040 | 11.5583 | 0.7961 | 110.340 |
| 15 | 5 | 4 | 1.9740 | 9.2647 | 0.7064 | 10.7099 |
| 16 | | 5 | 2.0290 | 9.7887 | 0.7328 | 21.9864 |
| 17 | | 6 | 2.0596 | 10.0855 | 0.7632 | 24.0855 |
| 18 | | 7 | 2.1128 | 10.6139 | 0.7675 | 40.0479 |
| 19 | | 8 | 2.1597 | 11.0903 | 0.7699 | 42.5799 |
| 20 | | 9 | 2.1943 | 11.4487 | 0.7763 | 72.4504 |
| 21 | | 10 | 2.2507 | 12.0438 | 0.7640 | 421.234 |
| 22 | 10 | 4 | 1.7637 | 9.1417 | 0.7159 | 169.622 |
| 23 | | 5 | 1.8094 | 9.6219 | 0.7455 | 86.3605 |
| 24 | | 6 | 1.8740 | 10.3211 | 0.7457 | 95.7238 |
| 25 | | 7 | 1.9214 | 10.8499 | 0.7508 | 159.780 |
| 26 | | 8 | 1.9391 | 11.0510 | 0.7726 | 330.072 |
| 27 | | 9 | 1.9664 | 11.3644 | 0.7820 | 185.579 |
| 28 | | 10 | 2.0018 | 11.7771 | 0.7813 | 232.406 |

Table 3. Egg packing results for test case 2

| Problem number | Number of container sides | Number of eggs | Objective | Area of optimized container | Packing fraction | Time (sec.) |
|---|---|---|---|---|---|---|
| 29 | 3 | 4 | 1.8381 | 4.3887 | 0.7457 | 15.1797 |
| 30 | | 5 | 1.9093 | 4.7357 | 0.7574 | 18.5818 |
| 31 | | 6 | 1.9336 | 4.8566 | 0.7924 | 23.9995 |
| 32 | | 7 | 1.9881 | 5.1343 | 0.7933 | 33.1757 |
| 33 | | 8 | 2.0180 | 5.2903 | 0.8070 | 47.1203 |
| 34 | | 9 | 2.0566 | 5.4945 | 0.8088 | 65.3648 |
| 35 | | 10 | 2.0913 | 5.6816 | 0.8098 | 75.5452 |
| 36 | 4 | 4 | 1.4386 | 4.1391 | 0.7906 | 17.1312 |
| 37 | | 5 | 1.5025 | 4.5148 | 0.7944 | 20.1131 |
| 38 | | 6 | 1.5461 | 4.7809 | 0.8050 | 24.5677 |
| 39 | | 7 | 1.5728 | 4.9476 | 0.8232 | 43.7241 |
| 40 | | 8 | 1.6144 | 5.2123 | 0.8191 | 47.2367 |
| 41 | | 9 | 1.6396 | 5.3765 | 0.8265 | 72.3506 |
| 42 | | 10 | 1.6713 | 5.5863 | 0.8236 | 99.1121 |
| 43 | 5 | 4 | 1.3064 | 4.0581 | 0.8064 | 15.3050 |
| 44 | | 5 | 1.3535 | 4.3560 | 0.8234 | 22.0519 |
| 45 | | 6 | 1.3909 | 4.6001 | 0.8366 | 29.7306 |
| 46 | | 7 | 1.4281 | 4.8492 | 0.8399 | 40.7732 |
| 47 | | 8 | 1.4614 | 5.0780 | 0.8407 | 52.5977 |
| 48 | | 9 | 1.4904 | 5.2815 | 0.8414 | 72.7178 |
| 49 | | 10 | 1.5118 | 5.4340 | 0.8467 | 85.8090 |
| 50 | 10 | 4 | 1.1624 | 3.9708 | 0.8241 | 53.9213 |
| 51 | | 5 | 1.1979 | 4.2170 | 0.8505 | 79.3888 |
| 52 | | 6 | 1.2422 | 4.5350 | 0.8486 | 105.983 |
| 53 | | 7 | 1.2777 | 4.7978 | 0.8489 | 133.923 |
| 54 | | 8 | 1.2974 | 4.9472 | 0.8630 | 177.798 |
| 55 | | 9 | 1.3255 | 5.1635 | 0.8606 | 224.357 |
| 56 | | 10 | 1.3481 | 5.3408 | 0.8615 | 294.655 |

Table 4. Egg packing results for test case 3

| Problem number | Number of container sides | Number of eggs | Objective | Area of optimized container | Packing fraction | Time (sec.) |
|---|---|---|---|---|---|---|
| 57 | 3 | 4 | 2.5905 | 8.7175 | 0.7548 | 17.5593 |
| 58 | | 5 | 2.7513 | 9.8335 | 0.7331 | 28.0943 |
| 59 | | 6 | 2.8194 | 10.3258 | 0.7490 | 40.0847 |
| 60 | | 7 | 2.8858 | 10.8185 | 0.7564 | 65.4501 |
| 61 | | 8 | 2.9339 | 11.1818 | 0.7670 | 80.5654 |
| 62 | | 9 | 2.9797 | 11.5340 | 0.7738 | 103.953 |
| 63 | | 10 | 3.0536 | 12.1128 | 0.7628 | 142.606 |
| 64 | 4 | 4 | 2.1050 | 8.8621 | 0.7425 | 21.2198 |
| 65 | | 5 | 2.1562 | 9.2988 | 0.7753 | 33.0706 |
| 66 | | 6 | 2.2376 | 10.0139 | 0.7723 | 41.2511 |
| 67 | | 7 | 2.2963 | 10.5457 | 0.7759 | 59.5441 |
| 68 | | 8 | 2.3385 | 10.9372 | 0.7841 | 82.6317 |
| 69 | | 9 | 2.4068 | 11.5853 | 0.7704 | 111.078 |
| 70 | | 10 | 2.4268 | 11.7785 | 0.7845 | 143.685 |
| 71 | 5 | 4 | 1.9439 | 8.9847 | 0.7324 | 22.6133 |
| 72 | | 5 | 1.9767 | 9.2901 | 0.7760 | 36.1682 |
| 73 | | 6 | 2.0239 | 9.7393 | 0.7941 | 49.2482 |
| 74 | | 7 | 2.0797 | 10.2841 | 0.7957 | 70.5491 |
| 75 | | 8 | 2.1322 | 10.8094 | 0.7934 | 95.7088 |
| 76 | | 9 | 2.1946 | 11.4508 | 0.7794 | 142.376 |
| 77 | | 10 | 2.2212 | 11.7306 | 0.7877 | 163.732 |
| 78 | 10 | 4 | 1.7507 | 9.0073 | 0.7305 | 84.8077 |
| 79 | | 5 | 1.7772 | 9.2819 | 0.7767 | 177.210 |
| 80 | | 6 | 1.8157 | 9.6893 | 0.7982 | 251.449 |
| 81 | | 7 | 1.8767 | 10.3510 | 0.7905 | 246.120 |
| 82 | | 8 | 1.8995 | 10.6040 | 0.8088 | 338.850 |
| 83 | | 9 | 1.9526 | 11.2051 | 0.7965 | 483.244 |
| 84 | | 10 | 1.9793 | 11.5138 | 0.8025 | 632.362 |

Table 5. Egg packing results for test case 4

| Problem number | Number of container sides | Number of eggs | Objective | Area of optimized container | Packing fraction | Time (sec.) |
|---|---|---|---|---|---|---|
| 85 | 3 | 4 | 1.8583 | 4.4860 | 0.7334 | 2.7068 |
| 86 | | 5 | 1.9052 | 4.7153 | 0.7645 | 4.1626 |
| 87 | | 6 | 1.9562 | 4.9708 | 0.7779 | 5.3980 |
| 88 | | 7 | 2.0110 | 5.2534 | 0.7788 | 8.0164 |
| 89 | | 8 | 2.0703 | 5.5678 | 0.7701 | 10.4752 |
| 90 | | 9 | 2.0871 | 5.6585 | 0.7887 | 14.5303 |
| 91 | | 10 | 2.1150 | 5.8110 | 0.7950 | 19.1139 |
| 92 | 4 | 4 | 1.4501 | 4.2056 | 0.7823 | 3.6087 |
| 93 | | 5 | 1.5141 | 4.5850 | 0.7862 | 5.8385 |
| 94 | | 6 | 1.5525 | 4.8208 | 0.8021 | 7.8616 |
| 95 | | 7 | 1.5782 | 4.9815 | 0.8213 | 9.8326 |
| 96 | | 8 | 1.6043 | 5.1478 | 0.8330 | 14.8706 |
| 97 | | 9 | 1.6495 | 5.4419 | 0.8201 | 19.4581 |
| 98 | | 10 | 1.6936 | 5.7365 | 0.8054 | 23.8477 |
| 99 | 5 | 4 | 1.3073 | 4.0632 | 0.8097 | 33.6934 |
| 100 | | 5 | 1.3527 | 4.3503 | 0.8286 | 52.4342 |
| 101 | | 6 | 1.3897 | 4.5919 | 0.8421 | 89.9398 |
| 102 | | 7 | 1.4305 | 4.8655 | 0.8409 | 107.572 |
| 103 | | 8 | 1.4617 | 5.0797 | 0.8441 | 153.185 |
| 104 | | 9 | 1.4905 | 5.2822 | 0.8449 | 919.116 |
| 105 | | 10 | 1.5286 | 5.5554 | 0.8316 | 259.104 |
| 106 | 10 | 4 | 1.1646 | 3.9860 | 0.8254 | 128.008 |
| 107 | | 5 | 1.2081 | 4.2896 | 0.8403 | 183.518 |
| 108 | | 6 | 1.2489 | 4.5842 | 0.8435 | 281.481 |
| 109 | | 7 | 1.2769 | 4.7915 | 0.8539 | 954.583 |
| 110 | | 8 | 1.3110 | 5.0512 | 0.8489 | 486.143 |
| 111 | | 9 | 1.3330 | 5.2218 | 0.8546 | 622.763 |
| 112 | | 10 | 1.3559 | 5.4027 | 0.8551 | 1435.69 |

Table 6. Egg packing results for test case 5

| Problem number | Number of container sides | Number of eggs | Objective | Area of optimized container | Packing fraction | Time (sec.) |
|---|---|---|---|---|---|---|
| 113 | 3 | 4 | 2.5323 | 8.3301 | 0.8026 | 769.011 |
| 114 | | 5 | 2.6901 | 9.4009 | 0.7785 | 888.957 |
| 115 | | 6 | 2.7724 | 9.9844 | 0.7857 | 1376.57 |
| 116 | | 7 | 2.8754 | 10.7403 | 0.7724 | 4189.08 |
| 117 | | 8 | 2.9253 | 11.1165 | 0.7817 | 3603.15 |
| 118 | | 9 | 3.0084 | 11.7566 | 0.7689 | 5123.83 |
| 119 | | 10 | 3.0549 | 12.1234 | 0.7717 | 6635.19 |
| 120 | 4 | 4 | 2.1002 | 8.8214 | 0.7579 | 806.650 |
| 121 | | 5 | 2.1920 | 9.6100 | 0.7615 | 408.624 |
| 122 | | 6 | 2.2683 | 10.2900 | 0.7623 | 2353.17 |
| 123 | | 7 | 2.2841 | 10.4342 | 0.7950 | 1886.77 |
| 124 | | 8 | 2.3317 | 10.8733 | 0.7992 | 5462.55 |
| 125 | | 9 | 2.3803 | 11.3316 | 0.7978 | 5119.08 |
| 126 | | 10 | 2.4475 | 11.9800 | 0.7809 | 7125.15 |
| 127 | 5 | 4 | 1.9254 | 8.8140 | 0.7586 | 753.826 |
| 128 | | 5 | 1.9641 | 9.1721 | 0.7979 | 1231.89 |
| 129 | | 6 | 2.0435 | 9.9285 | 0.7901 | 2060.84 |
| 130 | | 7 | 2.1054 | 10.5396 | 0.7871 | 3801.65 |
| 131 | | 8 | 2.1549 | 11.0405 | 0.7871 | 3512.30 |
| 132 | | 9 | 2.1745 | 11.2429 | 0.8040 | 11414.5 |
| 133 | | 10 | 2.2179 | 11.6956 | 0.7999 | 7243.70 |
| 134 | 10 | 4 | 1.7674 | 9.1800 | 0.7283 | 144.602 |
| 135 | | 5 | 1.7711 | 9.2187 | 0.7938 | 220.932 |
| 136 | | 6 | 1.8089 | 9.6169 | 0.8157 | 386.938 |
| 137 | | 7 | 1.8542 | 10.1038 | 0.8210 | 2241.72 |
| 138 | | 8 | 1.8674 | 10.2481 | 0.8479 | 1022.57 |
| 139 | | 9 | 1.9245 | 10.8852 | 0.8305 | 844.129 |
| 140 | | 10 | 1.9651 | 11.3488 | 0.8243 | 6043.67 |

Table 7. Egg packing results for test case 6

| Problem number | Number of container sides | Number of eggs | Objective | Area of optimized container | Packing fraction | Time (sec.) |
|---|---|---|---|---|---|---|
| 141 | 3 | 4 | 1.8521 | 4.4561 | 0.7502 | 19.9478 |
| 142 | 3 | 5 | 1.8946 | 4.6629 | 0.7847 | 61.6953 |
| 143 | 3 | 6 | 1.9510 | 4.9449 | 0.7932 | 226.774 |
| 144 | 3 | 7 | 2.0160 | 5.2798 | 0.7856 | 484.932 |
| 145 | 3 | 8 | 2.0764 | 5.6007 | 0.7758 | 808.216 |
| 146 | 3 | 9 | 2.0904 | 5.6764 | 0.7963 | 223.154 |
| 147 | 3 | 10 | 2.1299 | 5.8932 | 0.7937 | 1186.28 |
| 148 | 4 | 4 | 1.4592 | 4.2586 | 0.7850 | 1021.43 |
| 149 | 4 | 5 | 1.4992 | 4.4952 | 0.8140 | 1230.49 |
| 150 | 4 | 6 | 1.5558 | 4.8410 | 0.8102 | 1612.22 |
| 151 | 4 | 7 | 1.6041 | 5.1463 | 0.8060 | 2545.30 |
| 152 | 4 | 8 | 1.6322 | 5.3282 | 0.8154 | 3803.67 |
| 153 | 4 | 9 | 1.6629 | 5.5305 | 0.8173 | 5969.19 |
| 154 | 4 | 10 | 1.6882 | 5.6997 | 0.8207 | 6740.09 |
| 155 | 5 | 4 | 1.3348 | 4.2363 | 0.7891 | 2207.20 |
| 156 | 5 | 5 | 1.3775 | 4.5118 | 0.8110 | 3220.68 |
| 157 | 5 | 6 | 1.4176 | 4.7777 | 0.8210 | 4472.20 |
| 158 | 5 | 7 | 1.4462 | 4.9729 | 0.8341 | 6903.90 |
| 159 | 5 | 8 | 1.4881 | 5.2651 | 0.8252 | 6838.76 |
| 160 | 5 | 9 | 1.5067 | 5.3975 | 0.8374 | 20883.4 |
| 161 | 5 | 10 | 1.5293 | 5.5610 | 0.8411 | 13502.4 |
| 162 | 10 | 4 | 1.1859 | 4.1332 | 0.8088 | 192.373 |
| 163 | | 5 | 1.2239 | 4.4021 | 0.8312 | 1399.74 |
| 164 | | 6 | 1.2617 | 4.6783 | 0.8384 | 1937.92 |
| 165 | | 7 | 1.2917 | 4.9038 | 0.8458 | 6500.79 |
| 166 | | 8 | 1.3307 | 5.2038 | 0.8349 | 8730.52 |
| 167 | | 9 | 1.3436 | 5.3057 | 0.8519 | 3542.88 |
| 168 | | 10 | 1.3695 | 5.5118 | 0.8486 | 22030.7 |

Table 8. Egg packing results for test case 7

| Problem number | Number of container sides | Number of eggs | Objective | Area of optimized container | Packing fraction | Time (sec.) |
|---|---|---|---|---|---|---|
| 169 | 3 | 4 | 2.6259 | 8.9576 | 0.7324 | 4.4808 |
| 170 | | 5 | 2.8097 | 10.2554 | 0.7014 | 8.1396 |
| 171 | | 6 | 2.8555 | 10.5919 | 0.7289 | 14.4646 |
| 172 | | 7 | 2.9837 | 11.5647 | 0.7067 | 17.0902 |
| 173 | | 8 | 2.9682 | 11.4449 | 0.7488 | 44.3698 |
| 174 | | 9 | 3.0139 | 11.7999 | 0.7562 | 96.7891 |
| 175 | | 10 | 3.0565 | 12.1355 | 0.7615 | 825.447 |
| 176 | 4 | 4 | 2.1203 | 8.9910 | 0.7297 | 25.2150 |
| 177 | | 5 | 2.1505 | 9.2497 | 0.7776 | 42.6442 |
| 178 | | 6 | 2.2505 | 10.1294 | 0.7622 | 254.799 |
| 179 | | 7 | 2.3094 | 10.6664 | 0.7663 | 648.593 |
| 180 | | 8 | 2.3073 | 10.6475 | 0.8049 | 491.954 |
| 181 | | 9 | 2.3776 | 11.3063 | 0.7892 | 1249.69 |
| 182 | | 10 | 2.4305 | 11.8150 | 0.7821 | 2829.13 |
| 183 | 5 | 4 | 1.9532 | 9.0707 | 0.7233 | 5.1273 |
| 184 | | 5 | 2.0343 | 9.8391 | 0.7311 | 9.3015 |
| 185 | | 6 | 2.0398 | 9.8929 | 0.7804 | 17.5677 |
| 186 | | 7 | 2.1284 | 10.7706 | 0.7588 | 73.7508 |
| 187 | | 8 | 2.1273 | 10.7596 | 0.7965 | 43.6279 |
| 188 | | 9 | 2.1929 | 11.4336 | 0.7804 | 492.502 |
| 189 | | 10 | 2.2566 | 12.1079 | 0.7632 | 185.712 |
| 190 | 10 | 4 | 1.7538 | 9.0400 | 0.7257 | 75.2667 |
| 191 | | 5 | 1.7759 | 9.2686 | 0.7761 | 121.193 |
| 192 | | 6 | 1.8182 | 9.7154 | 0.7947 | 196.941 |
| 193 | | 7 | 1.8671 | 10.2456 | 0.7977 | 1287.00 |
| 194 | | 8 | 1.9274 | 10.9182 | 0.7849 | 1409.65 |
| 195 | | 9 | 1.9630 | 11.3250 | 0.7879 | 1761.68 |
| 196 | | 10 | 1.9846 | 11.5755 | 0.7983 | 6089.94 |

Table 9. Egg packing results for test case 8

| Problem number | Number of container sides | Number of eggs | Objective | Area of optimized container | Packing fraction | Time (sec.) |
|---|---|---|---|---|---|---|
| 197 | 3 | 4 | 2.8763 | 10.7470 | 0.7188 | 73.1147 |
| 198 |   | 5 | 2.9749 | 11.4967 | 0.7365 | 123.432 |
| 199 |   | 6 | 3.1845 | 13.1732 | 0.6897 | 151.958 |
| 200 |   | 7 | 3.2575 | 13.7846 | 0.6975 | 227.377 |
| 201 |   | 8 | 3.4068 | 15.0769 | 0.6685 | 298.352 |
| 202 |   | 9 | 3.5440 | 16.3158 | 0.6429 | 358.151 |
| 203 |   | 10 | 3.7169 | 17.9468 | 0.6052 | 669.925 |
| 204 | 4 | 4 | 2.2294 | 9.9408 | 0.7771 | 67.3030 |
| 205 |   | 5 | 2.3802 | 11.3307 | 0.7473 | 125.285 |
| 206 |   | 6 | 2.4112 | 11.6275 | 0.7813 | 163.078 |
| 207 |   | 7 | 2.4629 | 12.1321 | 0.7925 | 225.838 |
| 208 |   | 8 | 2.5937 | 13.4545 | 0.7491 | 311.393 |
| 209 |   | 9 | 2.7438 | 15.0572 | 0.6967 | 433.157 |
| 210 |   | 10 | 2.9695 | 17.6357 | 0.6159 | 510.207 |
| 211 | 5 | 4 | 2.1035 | 10.5206 | 0.7343 | 82.7059 |
| 212 |   | 5 | 2.2389 | 11.9184 | 0.7104 | 127.430 |
| 213 |   | 6 | 2.2976 | 12.5509 | 0.7238 | 189.892 |
| 214 |   | 7 | 2.3291 | 12.8980 | 0.7454 | 299.622 |
| 215 |   | 8 | 2.3904 | 13.5858 | 0.7418 | 312.729 |
| 216 |   | 9 | 2.5637 | 15.6268 | 0.6713 | 753.623 |
| 217 |   | 10 | 2.7022 | 17.3613 | 0.6256 | 757.450 |
| 218 | 10 | 4 | 1.9398 | 11.0588 | 0.6986 | 139.701 |
| 219 |   | 5 | 2.0385 | 12.2126 | 0.6933 | 210.233 |
| 220 |   | 6 | 2.0892 | 12.8277 | 0.7082 | 351.298 |
| 221 |   | 7 | 2.1804 | 13.9716 | 0.6882 | 437.045 |
| 222 |   | 8 | 2.3291 | 15.9423 | 0.6322 | 538.977 |
| 223 |   | 9 | 2.3294 | 15.9465 | 0.6578 | 888.381 |
| 224 |   | 10 | 2.4129 | 17.1113 | 0.6347 | 1492.95 |

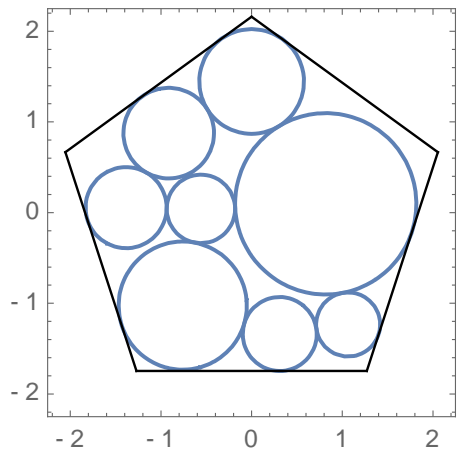

Test case 1

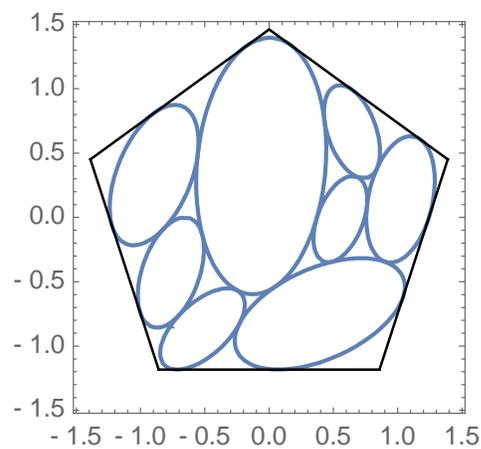

Test case 2

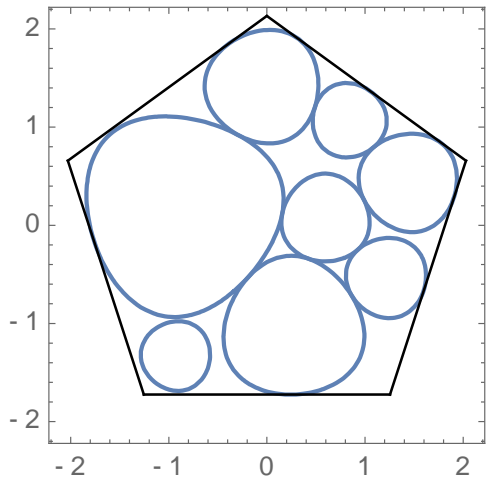
Test case 3

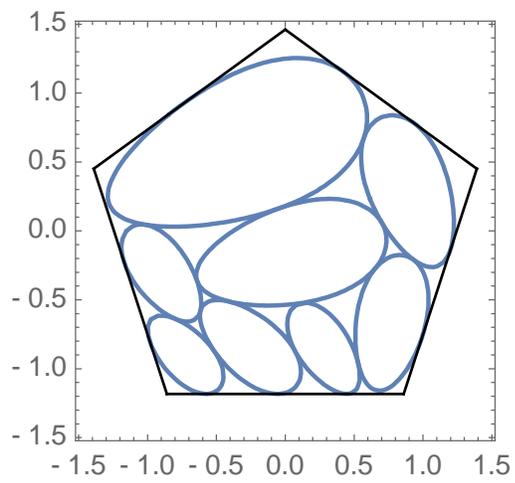
Test case 4

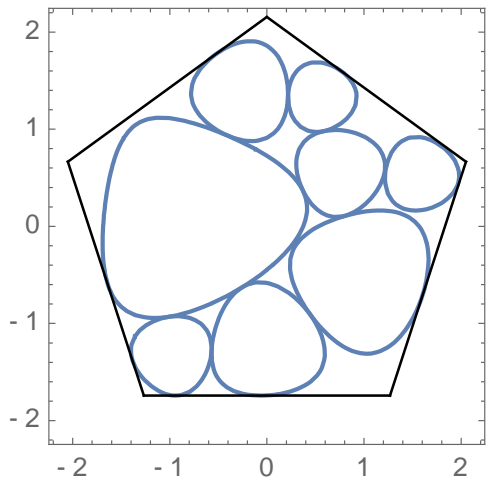
Test case 5

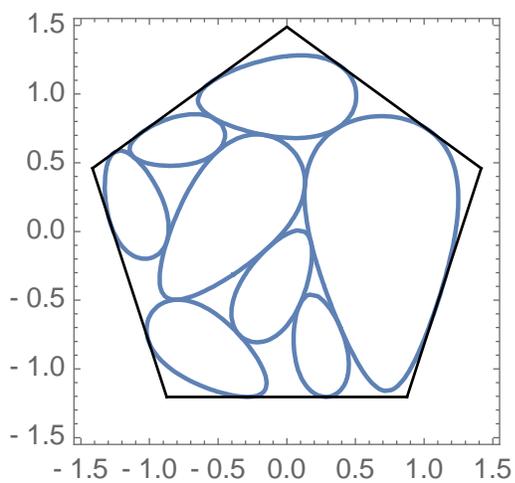
Test case 6

Figure 5. Packing configurations found for $n = 8$ and $m = 5$.

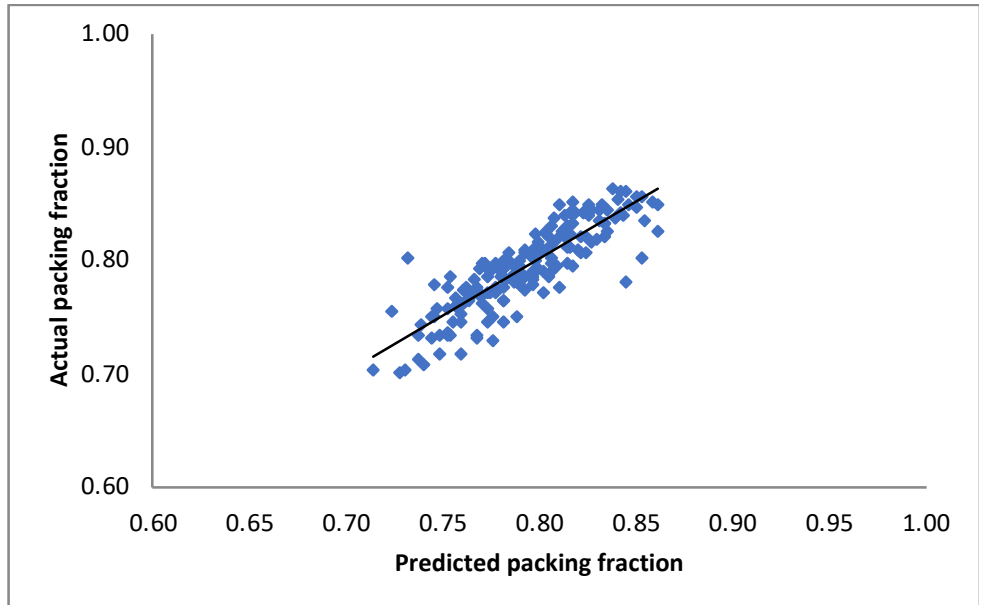

Figure 6. Actual (computed) vs. predicted packing fractions.

## 5 Conclusions

The efficient packing, arrangement, or configuration design of objects is required across a vast range of engineering and scientific applications. In this work, we present a model development approach to the challenging problem of packing convex planar objects. Specifically, we introduce the general problem-class of packing eggs (defined here as exponentially distorted ellipses) into optimized regular polygons. The numerical solution approach is based on the use of embedded Lagrange multipliers. We produce credible packings for all test problems considered, at the expense of an overall reasonable computational effort. Our embedded Lagrange multipliers based modeling approach is applicable to objects defined by *any* convex curve described by a function which is always positive for points outside the curve. This observation can lead to model development and solution strategies to handle further very general packing problems.

## References

Alder, B.J., Wainwright, T.E., 1957. Phase transition for a hard sphere system. *Journal of Chemical Physics* 27, 1208-1209.

Alt, H., 2016. Computational aspects of packing problems. In: *The Algorithmics Column*, *Bulletin of EATCS* 118. European Association for Theoretical Computer Science, www.eatcs.org.

Anjos, M.F. and Vieira, M.V.C., 2017. Mathematical optimization approaches for facility layout problems: The state-of-the-art and future research directions. *European Journal of Operational Research* 261, 1–16.

Bernal, J.D., 1959. Geometrical approach to the structure of liquids. *Nature* 183, 141-147.

Birgin, E.G., Lobato, R.D., Martínez, J.M., 2016. Packing ellipsoids by nonlinear optimization. *Journal of Global Optimization* 65, 709-743.

Birgin, E.G., Lobato, R.D., Martínez, J.M., 2017. A nonlinear programming model with implicit variables for packing ellipsoids. *Journal of Global Optimization* 68, 467-499.

Black, K., Chakrapani, C., Castillo, I., 2014. *Business Statistics for Contemporary Decision Making* (2nd Canadian Edition). Wiley & Sons Canada, Toronto, ON, Canada.

Bennell, J.A. and Oliveira, J.F., 2008. The geometry of nesting problems: A tutorial. *European Journal of Operational Research* 184, 397–415.

Bennell, J.A., Scheithauer, G., Stoyan, Y., and Romanova, T., 2010. Tools of mathematical modeling of arbitrary object packing problems. *Annals of Operations Research* 179, 343–368.

Castillo, I., Kampas, F.J., Pintér, J.D., 2008. Solving circle packing problems by global optimization: Numerical results and industrial applications. *European Journal of Operational Research* 191, 786-802.

Chaikin, P., 2000. Thermodynamics and hydrodynamics of hard spheres: the role of gravity. In: Cates, M.E., Evans, M.R., Eds. *Soft and Fragile Matter: Nonequilibrium Dynamics, Metastability and Flow*, Vol. 53. Institute of Physics Publishing.

Cheng, Z.D., Russell, W.B., Chaikin, P.M., 1999. Controlled growth of hard-sphere colloidal crystals. *Nature* 401, 893-895.

Chernov, N., Stoyan, Yu., Romanova, T., 2010. Mathematical model and efficient algorithms for object packing problem. *Computational Geometry* 43, 535–553.

Cohn, H., 2010. Order and disorder in energy minimization. *Proceedings of the International Congress of Mathematicians*, pp. 2416-2443. Hindustan Book Agency, New Delhi, India.

Conway, J.H., 1995. Sphere packings, lattices, codes, and greed. *Proceedings of the International Congress of Mathematicians*, pp. 45-55. Birkhäuser Verlag, Basel, Switzerland.

Dowsland, K.A. and Dowsland, W.B., 1992. Packing problems. *European Journal of Operational Research* 56, 2-14.

Edwards, S.F., 1994. The role of entropy in the specification of a powder. In: Mehta, A., Ed. *Granular Matter: An Interdisciplinary Approach*. Springer, New York.

Fasano, G. (2014) *Solving Non-standard Packing Problems by Global Optimization and Heuristics*. Springer International Publishing AG, Cham, Switzerland.

Fasano, G. and Pintér, J.D., Eds. (2015) *Optimized Packings with Applications.* Springer International Publishing AG, Cham, Switzerland.

Galiev, S.I., Lisafina, M.S., 2013. Numerical optimization methods for packing equal orthogonally oriented ellipses in a rectangular domain. *Computational Mathematics and Mathematical Physics* 53, 1748-1762.

Hifi, M., M'Hallah, R., 2009. A literature review on circle and sphere packing problems: models and methodologies. *Advances in Operations Research.*
DOI: 10.1155/2009/150624.

Ipopt (2018) https://projects.coin-or.org/Ipopt. The developers of Ipopt are listed at https://projects.coin-or.org/Ipopt/browser/trunk/Ipopt/AUTHORS.

Jaeger, H.M., Nagel, S.R., 1992. Physics of the granular state. *Science* 255, 1523-1531.

Jaeger, H.M., Nagel, S.R., Behringer, R.P., 1996. Granular solids, liquids, and gases. *Reviews of Modern Physics* 68, 1259-1273.

Jensen, F., 2017. *Introduction to Computational Chemistry.* (2$^{nd}$ Edition). John Wiley & Sons Ltd, Chichester, England.

Kallrath, J., 2017. Packing ellipsoids into volume-minimizing rectangular boxes. *Journal of Global Optimization* 67, 151-185.

Kallrath, J., Rebennack, S., 2014. Cutting ellipses from area-minimizing rectangles. *Journal of Global Optimization* 59, 405-437.

Kampas, F.J., Pintér, J.D., Castillo, I., 2017. Optimal packing of general ellipses in a circle. In: Takáč, M. and Terlaky, T., Eds., *Modeling and Optimization: Theory and Applications* (*MOPTA 2016 Proceedings*), pp. 23-38. Springer International Publishing AG, Cham, Switzerland.

Kampas, F.J., Castillo, I., Pintér, J.D., 2018. Optimized ellipse packings in regular polygons using embedded Lagrange multipliers. Research report available for download from *Optimization Online*; www.optimization-online.org/DB_FILE/2016/03/5348.pdf. (Submitted for publication)

Kellis, M., 2016. *Computational Biology: Genomes, Networks, Evolution.* An online textbook for MIT Course 6.047/6.878. Available at
https://ocw.mit.edu/ans7870/6/6.047/f15/MIT6_047F15_Compiled.pdf.

Kleijnen, J.P.C., 2015. *Design and Analysis of Simulation Experiments.* (2nd Edition). Springer US, New York.

Köller, J. (2018) Egg curves and ovals. http://www.mathematische-basteleien.de/eggcurves.htm.

Landau, R.H., Páez, M.J., and Bordeianu, C.C. 2012. *Computational Physics – Problem Solving with Computers.* Copyright © 2012 by Landau, Páez, and Bordeianu. Copyright © WILEY-VCH Verlag GmbH & Co. KGaA.

Lodi, A., Martello, S., and Vigo, D., 2002. Heuristic algorithms for the three-dimensional bin packing problem. *European Journal of Operational Research* 141, 410–420.

López, C.O. and Beasley, J.E., 2011. A heuristic for the circle packing problem with a variety of containers. *European Journal of Operational Research* 214, 512–525.

Newman, M., 2012. *Computational Physics.* CreateSpace Independent Publishing Platform.

O'Neil, S.T., 2017. *A Primer for Computational Biology*. Oregon State University Press, Corvallis, OR.

Pintér, J.D., Kampas, F.J., Castillo, I., 2017. Globally optimized packings of non-uniform size spheres in $R^d$: a computational study. *Optimization Letters*, DOI: 10.1007/s11590-017-1194-x.

Pisinger, D. and Sigurd, M., 2007. Using decomposition techniques and constraint programming for solving the two-dimensional bin-packing problem. *INFORMS Journal on Computing* 19 (1), 36–51.

Pusey, P.N., 1991. Colloidal suspensions. In: Hansen, J.P., Levesque, D., Zinnjustin, J., Eds. *Liquids, Freezing and Glass Transition*, Vol. 51 of Les Houches Summer School Session, 763-942. Elsevier Science Publishers, Amsterdam.

Rintoul, M.D., Torquato, S., 1996. Metastability and crystallization in hard-sphere systems. *Physical Review Letters* 77, 4198-4201.

Saunders, T.E. 2017. Imag(in)ing growth and form. *Mechanisms of Development* 145, 13–21.

Shannon, C.E., 1948. A mathematical theory of communication. *The Bell System Technical Journal 27,* 379–423 and 623–656 (July, October 1948).

Specht, E., 2018. http://www.packomania.com/.

Szabó, P.G., Markót, M.Cs, Csendes, T., Specht, E., Casado, L.G., García, I., 2007. *New Approaches to Circle Packing in a Square: With Program Codes.* Springer Science + Business Media, New York.

Thompson, D.W., 1917. *On Growth and Form.* Cambridge University Press.

Uhler, C., Wright, S.J., 2013. Packing ellipsoids with overlap. *SIAM Review* 55, 671-706.

Wäscher, G., Haußner, H., and Schumann, H., 2007. An improved typology of cutting and packing problems. *European Journal of Operational Research* 183, 1109–1130.

Wikipedia (2018) https://en.wikipedia.org/wiki/Oval.

Yamamoto, N. (2018) Equation of egg shaped curve.html. http://www.geocities.jp/nyjp07/index_egg_E.html